\documentclass[sn-basic]{sn-jnl}

\jyear{2022}%
\usepackage{amssymb,latexsym,amsmath,extarrows, mathrsfs, amsthm}
\usepackage{amssymb}
\usepackage{hyperref}
\usepackage{epstopdf}
\usepackage[utf8]{inputenc}
\theoremstyle{thmstyleone}%
\newtheorem{theorem}{Theorem}
\newtheorem{assumption}{Assumption}[section]

\theoremstyle{thmstyletwo}%

\theoremstyle{thmstylethree}%

\begin{document}

\title[Simultaneous Inference in Non-Sparse High-Dimensional Linear Models]{Simultaneous Inference  in Non-Sparse High-Dimensional Linear Models}


\author[1]{\fnm{Yanmei} \sur{Shi}}\email{2020020263@qdu.edu.cn}
\author[1]{\fnm{Zhiruo} \sur{Li}}\email{2020025413@qdu.edu.cn}
\author*[1]{\fnm{Qi} \sur{Zhang}}\email{qizhang@qdu.edu.cn}



\affil[1]{Institute of Mathematics and Statistics, Qingdao University, 308 Ningxia Road, Shinan District, Qingdao, Shandong, China}




\abstract{
  Inference and prediction under the sparsity assumption have been a hot research topic in recent years. However, in practice, the sparsity assumption is difficult to test, and more importantly  can usually be violated.
 In this paper, to study hypothesis test of any group of parameters under non-sparse high-dimensional linear models, we transform the null hypothesis to a testable moment condition and then use the self-normalization structure to construct moment test statistics under one-sample and two-sample cases, respectively.
 Compared to the one-sample case, the two-sample additionally requires a convolution condition.
 It is worth noticing that these test statistics contain  Modified Dantzig Selector, which simultaneously estimates model parameters and error variance without sparse assumption.
 Specifically, our method  can be extended to heavy tailed distributions of error for its robustness.
On very mild conditions, we show that the probability of  Type I error  is asymptotically equal to the nominal level $\alpha$ and the probability of Type II error is asymptotically 0. Numerical experiments indicate that our proposed method has good finite-sample performance.
}

\keywords{High-dimensional linear models, Non-sparse models, Group of parameters, Moment condition, Modified Dantzig Selector}

\pacs[MSC Classification]{62F03, 62F35, 62J15}


\maketitle

\section{Introduction}\label{sec1}

The rapid development of information technology has brought  important changes in data collection and data analysis.
 Nowadays, scientific researches in various fields increasingly rely on high-dimensional  observation data where the dimension $p$  is high, and the sample size $n$ is relatively small, that is $n\rightarrow \infty$ and $p/n\rightarrow \infty$.
 Modeling, inference and prediction for massive high-dimensional data have gradually become a research hotspot in statistics.
In addition, constrained by the current data scale and technology level, existing methods usually assume parameter sparsity, which means that the number of non-zero components of  parameter, denoted by $s$, is either fixed or increasing slowly as $s/n\rightarrow 0$ with $n \rightarrow \infty$\citep{2014Inference,2015Uniform,P2015High}.
Several penalized  minimization methods, including the Lasso \citep{Robert1996Regression}, Dantzig Selector \citep{Cand2007Rejoinder}, square root Lasso \citep{2011Square}, the ridge estimation \citep{B2013Statistical,2013Confidence} and scaled Lasso \citep{2013Confidence},  have been proposed and studied in sparse models.
In 2016, the Modified Dantzig Selector (MDS) was proposed by \cite{2016Linear}   to estimate model parameters under non-sparse conditions.
In general, penalization methods have been widely used  in many fields, for example biomedical imaging, disease tracking,  policy and market strategy.

Hypothesis test plays an important role in statistical theory and applications. Current progress in this field is also generally limited to the ultra-sparse case, i.e. $s=o(\sqrt{n}/logp)$.
\cite{2015Confidence} studied the confidence interval of high-dimensional linear regression with random design, and analyzed the expected length of the confidence interval and its adaptability at different sparsity levels;
 \cite{2014Confidence}  proposed a method for bias correction and showed that low-dimensional projection was an efficient way to construct
confidence intervals and hypothesis tests;
\cite{2015De} developed the de-biasing asymptotically Gaussian Lasso with the sparse level $\frac{n}{log^{2}p}$;
\cite{2014On} extended the estimation of one-component and low-dimensional components of parameter vectors to generalized linear models. \cite{2013A} established a two-sample global parameter homogeneity test which depended heavily on the model sparsity assumption. \cite{2016Testing} proposed the sLED test to compare the equivalence of two-sample covariance matrices;
\cite{Nicolas2018Two}  provided efficient asymptotic Type $\mathrm{I}$ error guarantees under sparse models.

It is pointed out that most existing inference methods sensitive to sparsity assumptions  might lose control of Type $\mathrm{I}$ error by \cite{2018Significance}, so they proposed approximately sparse CorrT which used reconstruct regression to transform the null hypothesis into a testable moment condition. They convolved two samples to translate the  homogeneity test of global parameter into a new moment condition \citep{2016Two}. They  further extended this moment method to hypothesis test of linear functionals of the regression parameters under a fully Gaussian design \citep{2016Linear}.

 The above literature  considered single parameters or global parameters,  to our knowledge, the current results on hypothesis testing for any group of non-sparse parameters are really rare.
  While in the field of genome-wide gene expression profiling research, \ $\beta \in \mathbb{R}^{p}$ is the whole genome, now  what we want to test is whether a group of parameters are equal to a given group of empirical parameters, that is  a group of genes are associated with a disease. For example, \cite{2016Mass}  pointed out that  genomes including B3GALNT1, C3orf62, TNFAIP1 and LTB gene linked to lung cancer with high probability. The test problem is as follows
\begin{align}
H_{0}: \beta_{G}=\beta_{G}^{0} \notag,
\end{align}
where  $\beta_{G} \in \beta$ is a group of parameters for  an arbitrary fixed  group $G\subseteq \{1,2,...,p\}$, and $\beta_{G}^{0}$ is a given parameter group.
Another example is to test whether two treatments $A$ and $B$ have the same effect on a disease.
Initially, we focused on global parameters, but the  treatment tends to affect a group  of cellular molecules levels.
Therefore, the problem degenerates into testing only the genetic cellular or molecular level $x_{ij,G}$ that is affected by the treatment, rather than the global level.
Supposing that $n_{1}$ patients receive treatment $A$ and $n_{2}$ patients receive treatment $B$,   $y_{i,j}$ denotes the response at the $x_{i,j}$ level for the $j$-th patient receiving treatment $i$ $(i=A \ or \  B)$, and the expectation of $y_{i,j}$ is $x_{i,j}^{T}\beta_{i}$ for $x_{i,j}\in \mathbb{R}^{p}$ and $\beta_{i}\in \mathbb{R}^{p}$.
The problem then turns into the following test
\begin{align}
H_{0}: \beta_{A,G}=\beta_{B,G} \notag,
\end{align}
where $\beta_{A,G}$ and $\beta_{B,G}$ denote the effect parameter groups of treatment A and treatment B at the genetic cellular or molecular level $x_{ij,G}$, respectively, for any fixed group $G\subseteq \{1,2,...,p\}$.

  Combining the methods in the above literature and the scenarios described,  we extend the reconstruction regression and convolutional regression methods to simultaneous inference of any group of parameters in the one- and two-sample cases.
 We transform the null hypothesis into moment conditions and construct  test statistics. The process does not rely heavily on accurate estimates of the parameters.
  Without imposing any sparsity assumptions, the test statistics guarantee asymptotic control of type I error and Type II error.

  The remainder of this paper is organized as follows: In Section \ref{section 2}, we develop a new method to perform one-sample simultaneous significance test, and  derive its size and power properties.
   Section \ref{section 3} extends the proposed results to the two-sample case, developing convolutional regression methods to perform two-sample homogeneity tests.
    In Section \ref{section 4}, Monte Carlo simulations are employed to demonstrate the excellent finite-sample performance of the proposed methods.
  The complete details of the theoretical proofs are contained in supplement.

\section{One-sample simultaneous inference} \label{section 2}
 In this section, we propose a new three-step method for  one-sample simultaneous significance test in non-sparse high-dimensional  models.  In the first step, we transform the null hypothesis into a moment condition by reconstructing regression.
  In the second step, the unknown parameters in the moment condition are estimated. In the third step, we normalize the moment condition to construct the test statistic, and derive its limiting distribution.
\subsection{Notations}
We first define the following notations. For a vector $v\in \mathbb{R}^{k}$, $\|v\|_{\infty}=\max\limits_{1\leq i \leq k}\vert v_{i}\vert$ and $\|v\|_{0}=\sum\limits_{i=1}^{k}\mathrm{I}(v_{i}\neq 0)$, where $\mathrm{I}( \cdot )$ denotes the indicator function. The  $v_{G}$ represents a sub-vector composed of $v_{j}, j\in G$ and $v_{-G}$ represents a sub-vector composed of $v_{j}, j \notin G$.
 For matrix $A$, its $(i,j)$ entry is denoted by $A_{i,j}$, and the $i$-th row and  $j$-th column are denoted by $a_{i}$ and $A_{j}$, respectively.
 The $A_{G}$ represents the submatrix composed of $A_{j}$, $ j\in G$ and $A_{-G}$ represents the submatrix composed of $A_{j}$, $ j\notin G$. $\|A\|_{\infty}=\max\vert A_{i,j} \vert$, where the maximum is taken over all $(i,j)$ indices.
 For two sequences $a_{n}, b_{n} >0$,  $a_{n}\asymp b_{n}$ means that there exist constants $C_{1}, C_{2} >0$ such that $\forall n$, $a_{n}\leq C_{1}b_{n}$ and $b_{n} \leq C_{2}a_{n}$. For two random quantities $X$ and $Y$ (scalars, vectors and matrices), $X\perp Y$ denotes the independence of $X$ and $Y$.

\subsection{Restructured statistical model}
Consider the following high-dimensional linear model
\begin{equation} \label{1.1}
y_{i}=x_{i}^{T}\beta+\varepsilon_{i}, i=1,2,...n, \tag{2.1}
\end{equation}
for a response $y_{i}\in\mathbb{R}$ and Gaussian design $x_{i}=(x_{i,G}^{T}, x_{i,-G}^{T})^{T}\in\mathbb{R}^{p}$, that is $x_{i}\sim N(0,\Sigma)$ with the unknown covariance matrix $\Sigma$.
 The noise vector $\varepsilon\in\mathbb{R}^{n}$ is internally uncorrelated and independent of the design matrix, satisfying $E(\varepsilon_{i})=0$, $E(\varepsilon_{i}^{2})=\sigma_{\varepsilon}^{2}$  with  $0<\sigma_{\varepsilon}^{2}<\infty$. The parameter vector $\beta=(\beta_{G}^{T},\beta_{-G}^{T})^{T}\in\mathbb{R}^{p}$ is unknown and is  allowed  $p\gg n$. We intend to test
\begin{equation}\label{1.2}
H_{0,G}: \beta_{G}=\beta_{G}^{0}  \tag{2.2},
\end{equation}
where $\beta_{G}$ is a parameter group we are interested in,  for the subvector $G\subseteq \{1,2,...,p\}$, and $\beta_{G}^{0}$ is given.

In high-dimensional dense models, it is difficult to directly construct consistent estimators of the parameters and determine their asymptotic distribution. So we now take a new approach to  test \eqref{1.2}.
First,  we rewrite the  model \eqref{1.1} into the following form
\begin{equation}
y=z^{T}\gamma+w^{T}\theta+\varepsilon  \notag,
\end{equation}
where $z=x_{G}\in \mathbb{R}^{k}$, $w=x_{-G}\in \mathbb{R}^{p-k}$,  $\gamma=\beta_{G}\in \mathbb{R}^{k}$ and $\theta=\beta_{-G}\in \mathbb{R}^{p-k}$. Then,
we introduce a pseudo-response $v:=y-z^{T}\beta_{G}^{0} $ to build the $reconstructed \ regression \  linear \ model$ as
  \begin{equation} \label{2.1}
 v=w^{T}\theta+e \tag{2.3},
  \end{equation}
   where the pseudo-error $e=z^{T}(\gamma-\beta_{G}^{0})+\varepsilon$ with $E(e)=z^{T}(\gamma-\beta_{G}^{0})$ and $\sigma_{e}^{2}=E(e^{2})=\sigma_{\varepsilon}^{2}$.  Obviously, $w$ and  $e$ are not correlated  under the null hypothesis $H_{0}$  in this new  model, and $e=\varepsilon$. However, under the alternative hypothesis $H_{1}$,  $e$ and  $w$ may be linear dependent through $z$, which is caused by the confounding effects of $w$ and $z$.
    So we consider the following linear model
   \begin{equation}\label{2.2}
   z=\pi^{T}w+u  \tag{2.4}.
   \end{equation}
 The $j$-th column $\pi_{j}\in\mathbb{R}^{p-k}$ of $\pi \in\mathbb{R}^{(p-k)\times k}$ is an unknown regression coefficient vector, and $u\in \mathbb{R}^{k}$ is internally uncorrelated,  satisfying $E(u)=0$ and $ E(uw^{T})=0$. In particular, when $u$ follows a Gaussian distribution,  $u\perp w$. In this article, we will assume that $\pi$ is sparse, in order to decouple the dependence between $z$ and $w$. In fact, sparse $\pi$ is a generalization of the sparsity condition on the precision matrix $\Omega_{X}=\Sigma_{X}^{-1}$, which is a typical regularity condition.

Considering the correlation between the pseudo-error $e$ in the new model \eqref{2.1} and the error $u$ in \eqref{2.2}, we construct the following moment condition
\begin{equation}
E(ue)
=E[uz^{T}(\gamma-\beta_{G}^{0})+u\varepsilon]
=E(uu^{T})(\gamma-\beta_{G}^{0}) \notag,
\end{equation}
Under the null hypothesis in \eqref{1.2}, $E(ue)=0$; on the contrary, $E(ue)\neq0$. So the hypothesis $H_{0,G}$ is equivalent to
\begin{equation} \label{2.3}
H_{0}:E[(z-\pi^{T}w)(v-w^{T}\theta)]=0
\tag{2.5}.
\end{equation}
 Since $\pi$ is sparse, we can get its consistent estimator $\hat{\pi}$.
  While it is difficult to obtain consistent estimator of  $\theta$ due to its non-sparseness.
  For any estimator  $\check{\theta}$ of $\theta$, we have
 \begin{align}
E(z-\hat{\pi}^{T}w)(v-w^{T}\check{\theta})
\rightarrow E(z-\pi^{T}w)(v-w^{T}\check{\theta})
= 0. \notag
 \end{align}
 Therefore, the above inner product structure alleviates the reliance on a good estimator of $\theta$.
 The moment condition in \eqref{2.3} contains the sparse parameter $\pi$ and the dense parameter $\theta$, which we want to replace with estimators.

\subsection{Modified Dantzig Selector} \label{section 2.2}
   Given the sample $X=(x_{1},x_{2},...,x_{n})^{T}\in\mathbb{R}^{n\times p}$, we define the new response vector  $V=Y-Z\beta_{G}^{0} \in \mathbb{R}^{n}$, new design matrices  $Z=X_{G}\in \mathbb{R}^{n\times k}$, $W=X_{-G}\in\mathbb{R}^{n\times(p-k)}$,  and $e=Z(\gamma-\beta_{G}^{0})+\varepsilon=(e_{1},e_{2},...e_{n})^{T}$.
Because this paper considers a non-sparse high-dimensional model, we use the MDS \citep{2016Linear}  to estimate   the unknown parameters and error variance $\sigma^{2}$   simultaneously.\

 The MDS estimator of $\theta$ is defined as follows
 \begin{equation}\tag{2.6} \label{theta_tilde}
\begin{split}
\tilde{\theta}=&\arg\min\limits_{\theta\in\mathbb{R}^{p-k}} \|\theta\|_{1} \notag
\\
s.t.& \ \|W^{T}(V-W\theta)\|_{\infty}\leq \eta\rho_{1}\sqrt{n}\|V\|_{2}
\notag \\
&V^{T}(V-W\theta)\geq \rho_{0}\rho_{1}\|V\|_{2}^{2}/2 \notag \\
&\rho_{1} \in [\rho_{0}, 1] \notag,
\end{split}
\end{equation}
where  $\rho_{1}=\sigma_{e}/ \sqrt{E(v)^{2}}$ and $\rho_{0} \in (0,1)$ is a lower bound for this ratio,
$\eta\asymp \sqrt{n^{-1}log p}$ and $\rho_{0}$  are tuning parameters.

Similarly,   $\tilde{\pi}_{j}\in \mathbb{R}^{p-k}$ for $j\in G$ is
\begin{equation} \tag{2.7} \label{pij_tilde}
\begin{split}
\tilde{\pi}_{j}=&\arg\min\limits_{\pi_{j}\in\mathbb{R}^{p}} \|\pi_{j}\|_{1}
\notag \\
s.t.& \|W^{T}(Z_{j}-W\pi_{j})\|_{\infty}\leq\eta\rho_{2}\sqrt{n}\|Z_{j}\|_{2}
\notag \\
&Z_{j}^{T}(Z_{j}-W\pi_{j})\geq \rho_{0}\rho_{2}\|Z_{j}\|_{2}^{2}/2 \notag \\
& \rho_{2}\in[\rho_{0},1] \notag,
\end{split}
\end{equation}
where $\rho_{2}=\sigma_{u}/\sqrt{E(z_{j})^{2}}$.
 Then the MDS estimator of $\pi$ is
\begin{equation} \label{pi_tilde}
\tilde{\pi}=(\tilde{\pi}_{j},j\in G)\in\mathbb{R}^{(p-k)\times k}.
\tag{2.8}
\end{equation}

\subsection{Simultaneous test} \label{ section 2.3}
 In this subsection, we construct a test statistic  using $\tilde{\theta}$ and $\tilde{\pi}$  to test the moment condition in \eqref{2.3} and take a  simulation method to obtain critical values of the rejection.

 By plugging in the estimators, the test statistic is
 \begin{equation}\label{jianyantongjiliang}
T_{n}=n^{-\frac{1}{2}}\hat{\sigma}_{e}^{-1}\|(Z-W\tilde{\pi})^{T}(V-W\tilde{\theta})\|_{\infty}
\tag{2.9},
\end{equation}
where  $\hat{\sigma}_{e}=\|V-W\tilde{\theta}\|_{2}/\sqrt{n}$.
 When the null hypothesis holds, the value of $T_{n}$ tends to be moderate; otherwise, the value of $T_{n}$ is large. Therefore, if $T_{n}$ is "too large",  the null hypothesis $H_{0}$ should be rejected and the alternative hypothesis $H_{1}$ should be accepted.
Due to the complex dependencies between  different terms of $(Z-W\tilde{\pi})^{T}(V-W\tilde{\theta})$ and high-dimensional non-sparse assumption, it is difficult to obtain the exact distribution of $T_{n}$.
However,  motivated by  Gaussian approximation method \citep{2014On}, we take the critical value  as a pre-specified quantile of the $l_{\infty}$ norm of a Gaussian vector with a zero-mean and a known covariance, which is easy to calculate by simulation.

Obviously, there is the following decomposition under the null hypothesis.
\begin{equation}
n^{-\frac{1}{2}}(Z-W\tilde{\pi})^{T}(V-W\tilde{\theta}) =
\Delta+n^{-\frac{1}{2}}U^{T}\hat{e} \notag,
\end{equation}
where $\hat{e}=V-W\tilde{\theta}$,
\begin{equation}
\Delta=n^{-\frac{1}{2}}(\pi-\tilde{\pi})^{T}W^{T}\hat{e} \notag,
 \end{equation}
  and $U\in \mathbb{R}^{n\times k}$ with
$ U_{j}=Z_{j}-W\pi_{j} $.
 We could show that $\|\Delta\|_{\infty}\hat{\sigma}_{e}^{-1}=o_{p}(1)$, see Lemma 2 in the supplement for details.
So the statistical properties of the test statistic $T_{n}$ is determined by $\|n^{-\frac{1}{2}}U^{T}\hat{e}\hat{\sigma}_{e}^{-1}\|_{\infty}=\|n^{-\frac{1}{2}}\sum\limits_{i=1}^{n}u_{i}\hat{e}_{i}\hat{\sigma}_{e}^{-1}\|_{\infty}$.

Under the null hypothesis, our constructed $U$ is independent of $(V, W)$, while  the MDS estimator $\tilde{\theta}$  is completely dependent on $(V,W)$, so $\hat{e}=V-W\tilde{\theta}$ is also only related to $(V,W)$. Therefore, $U$ and $\hat{e}$ are independent. Assuming that $U$ obeys a Gaussian distribution, that is $U\sim N(0,E(u_{1}u_{1}^{T}))$, where $u_{i}^{T}$ is the $i$-th row of $U$. Then
\begin{align}
E[n^{-\frac{1}{2}}\sum\limits_{i=1}^{n}u_{i}\hat{e}_{i}\hat{\sigma}_{e}^{-1}]=&E(\frac{1}{\|\hat{e}\|_{2}}\sum\limits_{i=1}^{n}u_{i}\hat{e}_{i})=0 \notag, \\
Var[n^{-\frac{1}{2}}\sum\limits_{i=1}^{n}u_{i}\hat{e}_{i}\hat{\sigma}_{e}^{-1}]
=&n^{-1}\hat{\sigma}_{e}^{-2}\sum\limits_{i=1}^{n}Var(u_{i}\hat{e}_{i})
=E(u_{1}u_{1}^{T}) \notag.
\end{align}
Therefore,
\begin{equation}
n^{-\frac{1}{2}}\hat{\sigma}_{e}^{-1}U^{T}\hat{e}\sim N(0, E(u_{1}u_{1}^{T}) ).
\notag
\end{equation}

Let $Q=E(u_{1}u_{1}^{T})$, then $n^{-\frac{1}{2}}\hat{\sigma}_{e}^{-1}U^{T}\hat{e}\sim N(0, Q ) $, and $Q$ is unknown. Consider $\hat{Q}=\frac{1}{n}\sum\limits_{i=1}^{n}\hat{u}_{i}\hat{u}_{i}^{T}$,
where $\hat{u}_{i}^{T}$  is the $i$-th row of $\hat{U}=Z-W\tilde{\pi} \in \mathbb{R}^{n\times k}$.
We introduce the function $\Gamma(x;A):=p(\|\xi\|_{\infty}\leq x) $, where $\xi\sim N(0, A)$.
The following Theorem \ref{theorem 2.1} shows that  $\Gamma(x; \hat{Q})$ can  asymptotically approximate the distribution of  $T_{n}$.
Moreover, $\Gamma(\cdot ; \hat{Q})$ can be easily simulated.
We summarize the above process into Algorithm \ref{1}.
\begin{algorithm}
\caption{One-sample simultaneous inference.} 
\label{1}
\begin{algorithmic}[1] 
\Require

Sample $(X, Y)$ and parameter $\beta_{G}^{0}$ with any subgroup $G \subseteq \{1,2,...,p \}$;
The nominal level $\alpha\in(0, 1)$ of the test.
\Ensure
Determine whether to reject the null hypothesis $H_{0}: \beta_{G}=\beta_{G}^{0}$.
\State Construct $Z=X_{G} $, $ W=X_{-G}$, $\gamma=\beta_{G}$ and $\theta=\beta_{-G}$.
The regression model is
$Y=Z\gamma+W\theta+\varepsilon$.
\State  Construct $V=Y-Z\beta_{G}^{0}$ and  $e=Z(\gamma-\beta_{G}^{0})+\varepsilon$. \\
The reconstructed model is $V=W\theta+e$.
\State  Determine the linear correlation model between $Z_{j}$ and $W$: $Z=W\pi+U$.
\State  Calculate $\tilde{\theta}$, $\tilde{\pi}$ and $\hat{Q}$, take $\eta\asymp\sqrt{n^{-1}logp}$.
\State  Calculate the test statistic $T_{n}=n^{-\frac{1}{2}}\hat{\sigma}_{e}^{-1}\|(Z-W\tilde{\pi})^{T}(V-W\tilde{\theta})\|_{\infty}$.
\State  Approximately calculate $\Gamma^{-1}(1-\alpha; \hat{Q})$.\\
\Return Reject the null hypothesis $H_{0}$ if and only if $T_{n}>\Gamma^{-1}(1-\alpha; \hat{Q}) $.
\end{algorithmic}
\end{algorithm}
\subsection{Theoretical results}
We introduce the size and power properties of the test statistic $T_{n}$ in this subsection while  imposing extremely weak assumptions when both $n$ and $p$ tend to $\infty$. The hypothesis testing problem is

\begin{align} \label{2.15}
H_{0}: \beta_{G}=\beta_{G}^{0} \ v.s. \ H_{1}: \beta_{G}=\beta_{G}^{0}+\mathbf{h} \tag{2.10},
\end{align}
where $\mathbf{h}\in \mathbb{R}^{k}$ is a nonzero vector. We first make the following assumptions.
\begin{assumption} \label{assumption 2.1}
Consider the model \eqref{1.1}. Suppose that the following hold: \\
(i) the design matrix $X$ follows a Gaussian distribution, that is $X\sim N(0,\Sigma)$;\\
(ii)  there exist constants $c$, $d$ $\in (0,+\infty)$ such that the eigenvalues of covariance matrix $\Sigma$ lie in $(c,d)$;\\
(iii)  the $\pi$ is sparse, which means $s_{\pi}=o(\sqrt{n/log^{3}p})$, where
$s_{\pi}=\max\limits_{j \in G}{\|\pi_{j}\|_{0}}$; \\
(iv) $s_{\pi}\|\beta_{G}-\beta_{G}^{0}\|_{0}+\|\theta\|_{0}=o(\sqrt{n}/logp)$; \\
(iv') $\|\beta_{G}-\beta_{G}^{0}\|_{0}=o_{p}(\sqrt{logp})$ and $\|\theta\|_{0}=o(\sqrt{n}/logp)$. \\
(v) there exist constants $\delta \  and\  \kappa_{1} \in (0,+\infty)$ such that: $E\vert \varepsilon \vert^{2+\delta}<\kappa_{1}$;\\
(vi) there exist constants $K_{1} \ and \ K_{2} >0$, which only depends on the constants $\delta \  and\  \kappa_{1} $, satisfying that
 $ \|\Sigma_{u}(\beta_{G}-\beta_{G}^{0})\|_{\infty}\geq\sqrt{n^{-1}logp}(K_{1}\|\beta_{G}-\beta_{G}^{0}\|_{2}+K_{2}) $.
\end{assumption}
Assumption  \ref{assumption 2.1} is relatively mild.
 Assumption  \ref{assumption 2.1}(i) is imposed to simplify the proofs.
 Assumption  \ref{assumption 2.1}(ii) is very standard in high-dimensional literature (see \cite{2014On}).
 Assumption \ref{assumption 2.1}(iii) imposes a sparsity condition on the regression coefficient vector $\pi_{j}$, rather than on $\beta$ and $\Sigma$ in the original model \eqref{1.1}.
 In Assumption \ref{assumption 2.1}(iv), a certain sparse structure is used to guarantee the asymptotic power of high-dimensional tests. But non-sparse $\beta$ is allowed.
For example, $\beta=(\gamma, \theta)^{T}$ is a dense parameter with $\|\gamma\|_{0}=o(\sqrt{log p})$ and $\|\theta\|_{0}=o(\sqrt{n}/log p)$, satisfying $\max \|\gamma-\beta_{G}^{0}\|_{0}=o(\sqrt{log p})$, then $ s_{\pi}\|\gamma-\beta_{G}^{0}\|_{0}+\|\theta\|_{0}=o(\sqrt{n}/logp)$. For sparse vector $\beta$, the sparsity condition $o(\sqrt{n}/logp)$ is consistent with the traditional test \citep{2013Two,2014On}.
 Note that Assumption \ref{assumption 2.1}(iv')  is   stricter than Assumption \ref{assumption 2.1}(iv), and relaxing it gives Assumption \ref{assumption 2.1}(iv).
 Assumption \ref{assumption 2.1}(v) is a regular moment condition.
 For the convenience of  proof, we impose a $l_{\infty}$ norm lower bound on the product of $\Sigma_{u}$ and the deviation in Assumption \ref{assumption 2.1}(vi).
Then  we provide the following results for  $T_{n}$.
\begin{theorem} \label{theorem 2.1}
 Let Assumption \ref{assumption 2.1}(i)-(iii) hold, when $n, p\rightarrow \infty$ with $logp=o(\sqrt{n})$, then under null hypothesis,
 \begin{align}
  P(T_{n}>\Gamma^{-1}(1-\alpha;\hat{Q}))\rightarrow\alpha, \forall\alpha\in(0,1)
 \tag{2.11},
 \end{align}
 where $\Gamma^{-1}(1-\alpha;\hat{Q})$  is the $1-\alpha$ quantile of the $\Gamma(x; \hat{Q})$.
\end{theorem}
Theorem \ref{theorem 2.1} gives the asymptotic distribution of $T_{n}$.
In contrast to existing methods, we do not make any sparsity assumptions about the model parameter $\beta$; other than that, there are no constraints on the distribution of the error $\varepsilon$ in model \eqref{1.1}.
The inner product  $n^{-\frac{1}{2}}(Z-W\tilde{\pi})^{T}(V-W\tilde{\theta})\hat{\sigma}_{e}^{-1}$  is normally distributed under $H_{0}$.
 In addition, we consider the power property of  $T_{n}$.

\begin{theorem} \label{theorem 2.2}
  Let Assumption \ref{assumption 2.1} holds, under the alternative hypothesis in \eqref{2.15}, when $n, p\rightarrow \infty$, with $logp=o(\sqrt{n})$, then the test in Algorithm \ref{1} is asymptotically  powerful, $i.e.$,
 \begin{align}
 P(T_{n}>\Gamma^{-1}(1-\alpha;\hat{Q}))\rightarrow 1, \notag \ {\forall} \alpha \in (0,1).
 \end{align}
 \end{theorem}
 Theorem \ref{theorem 2.2}  shows that no matter the distribution of model error,   the test of  any parameter group is asymptotically powerful, that is, the probability of Type II error  tends to 0.
\section{Two-sample partial equivalence inference}  \label{section 3}
\subsection{Partial equivalence inference}
 Here, we highlight the extension of the proposed method to  two-sample models. We consider the following two-sample linear regression models
 \begin{align}
y_{A}=x_{A}^{T}\beta_{A}+\varepsilon_{A} \tag{3.1}\label{1.3} \\
y_{B}=x_{B}^{T}\beta_{B}+\varepsilon_{B} \tag{3.2}\label{1.4},
\end{align}
for unknown parameters $\beta_{A}=(\beta_{A,G}^{T},\beta_{A,-G}^{T})^{T}\in\mathbb{R}^{p}$ and $\beta_{B}=(\beta_{B,G}^{T},\beta_{B,-G}^{T})^{T}\in\mathbb{R}^{p}$. We also consider Gaussian designs, that is   $x_{A}=(x_{A,G}^{T},x_{A,-G}^{T})^{T}\sim N(0,\Sigma_{A})$ and $x_{B}=(x_{B,G}^{T},x_{B,-G}^{T})^{T}\sim N(0,\Sigma_{B})$ with unknown covariance matrices $\Sigma_{A}$ and $\Sigma_{B}$. Furthermore, the errors $\varepsilon_{A}$ and $\varepsilon_{B}$ satisfy $E(\varepsilon_{A})=0$ and $E(\varepsilon_{A}^{2})=\sigma_{\varepsilon, A}^{2}$, $E(\varepsilon_{B})=0$ and $E(\varepsilon_{B}^{2})=\sigma_{\varepsilon, B}^{2}$ with  unknown error variances $\sigma_{\varepsilon, A}^{2}$ and $\sigma_{\varepsilon, B}^{2}$.
Our goal is to test whether partial parameters are  identical in \eqref{1.3}  and \eqref{1.4}, that is,
\begin{equation} \label{1.5}
H_{0}: \beta_{A,G}=\beta_{B,G} \tag{3.3},
\end{equation}
where $\beta_{A,G}$ and $\beta_{B,G}$ are  parameter groups we are interested in,  for the given subvector $G\subseteq \{1,2,...,p\}$.

Now we extend reconstruction and  convolutional regression methods to  the partial equivalence test.
Firstly, we reconstruct the models \eqref{1.3} and \eqref{1.4} to
\begin{align}
y_{A}=x_{A,G}^{T}\beta_{A,G}+x_{A,-G}^{T}\beta_{A,-G}+\varepsilon_{A} \tag{3.4} \label{yangben1chonggou} \\
y_{B}=x_{B,G}^{T}\beta_{B,G}+x_{B,-G}^{T}\beta_{B,-G}+\varepsilon_{B} \tag{3.5} \label{yangben2chonggou}.
\end{align}
Secondly, we convolve the variables of new models \eqref{yangben1chonggou} and \eqref{yangben2chonggou}.
 Define the convolution response $y=y_{A}+y_{B}$, the error $\varepsilon_{*}=\varepsilon_{A}+\varepsilon_{B}$, new design matrices $z=x_{A,G}-x_{B,G}\in \mathbb{R}^{k}$, $m=x_{A,G}+x_{B,G}\in \mathbb{R}^{k}$ and $w=(m,x_{A,-G},x_{B,-G})\in \mathbb{R}^{2p-k}$.
The $convolution \ regression \  model$ is as follows
\begin{equation} \label{3.3}
y=z^{T}\gamma_{*}+w^{T}\theta_{*}+\varepsilon_{*}, \tag{3.6}
\end{equation}
with unknown parameters $\gamma_{*}=(\beta_{A,G}-\beta_{B,G})/2\in \mathbb{R}^{k}$ and $\theta_{*}=[(\beta_{A,G}^{T}+\beta_{B,G}^{T})/2, \beta_{A,-G}, \beta_{B,-G}]^{T}\in \mathbb{R}^{2p-k}$. It can be seen that $\beta_{A,G}-\beta_{B,G}$ appears as regression coefficients in \eqref{3.3}.  Therefore, the null hypothesis in \eqref{1.5} is equivalent to
\begin{equation} \label{3.4}
H_{0}:\gamma_{*}=0 \tag{3.7}.
\end{equation}
Next, we construct the test statistic in the same way as one-sample case and estimate its parameters using the MDS method. The above process is summarized in Algorithm \ref{2}.

\begin{algorithm}
\caption{Two-sample partial homogeneity test.} 
\label{2}
\begin{algorithmic}[1] 
\Require
Sample $X_{A}=(x_{A,1},x_{A,2},...,x_{A,n})^{T}\in\mathbb{R}^{n\times p}$, $X_{B}=(x_{B,1},x_{B,2},...,x_{B,n})^{T}\in\mathbb{R}^{n\times p}$,
  $Y_{A}=(y_{A,1},y_{A,2},...,y_{A,n})^{T}\in\mathbb{R}^{n}$,
  and $Y_{B}=(y_{B,1},y_{B,2},...,y_{B,n})^{T}\in\mathbb{R}^{n}$;
the nominal level $\alpha\in(0, 1)$ of the test.
\Ensure
Determine whether to reject the null hypothesis $H_{0}: \beta_{A, G}=\beta_{B, G}$.
\State  Construct $Z=X_{A,G}-X_{B,G}$, $W=(X_{A,G}+X_{B,G},X_{A,-G},X_{B,-G})$ and $Y=Y_{A}+Y_{B}\in \mathbb{R}^{n}$. \\
The convolution regression model is
 $Y=Z\gamma_{*}+W\theta_{*}+\varepsilon_{*}$.
\State  Determine the linear correlation model between $Z$ and $W$: $Z=W\pi_{*}+U_{*}$.
\State  Calculate estimators $\tilde{\theta}_{*}$ and $\tilde{\pi}_{*}$  using MDS method and define $\hat{D}=n^{-1}\sum\limits_{i=1}^{n} \hat{u}_{*,i}\hat{u}_{*,i}^{T}$.
\State  Calculate the test statistic $S_{n}=n^{-\frac{1}{2}}\hat{\sigma}_{\varepsilon_{*}}^{-1}\|(Z-W\tilde{\pi}_{*})^{T}(Y-W\tilde{\theta}_{*})\|_{\infty} $,
 where $\hat{\sigma}_{\varepsilon_{*}}=\|Y-W\tilde{\theta}_{*}\|_{2}/\sqrt{n}$.
 \State  Approximately calculate $\Gamma^{-1}(1-\alpha; \hat{D})$ (by simulation).\\
\Return Reject the null hypothesis $H_{0}$ \eqref{3.4} if and only if $S_{n}>\Gamma^{-1}(1-\alpha; \hat{D}) $.
\end{algorithmic}
\end{algorithm}
\subsection{Theoretical results}
We now study the theoretical results of the test statistic $S_{n}$ under mild assumptions. The hypothesis testing problem is
\begin{align}
H_{0}: \beta_{A,G}=\beta_{B,G} \ v.s. H_{1}: \beta_{A,G}=\beta_{B,G}+\mathbf{h}, \tag{3.8} \label{result hypothesis}
\end{align}
where $\mathbf{h}\in \mathbb{R}^{k}$ is a nonzero vector. Similar to the one-sample case, we make the following assumptions.

\begin{assumption} \label{assumption 3.1}
Consider models \eqref{1.3} and \eqref{1.4}. Suppose that the following hold:\\
(i) the design vectors follow  Gaussian distribution,
 $x_{A}\sim N(0, \Sigma_{A})$, $x_{B}\sim N(0, \Sigma_{B})$; \\
 (ii) there exist constants $\kappa_{1}$, $\kappa_{2}$ such that the eigenvalues of \  $\Sigma_{A}$ and $\Sigma_{B}$ lie in $(\kappa_{1}, \kappa_{2})$; \\
 (iii) the $\pi_{*}$ is sparse, which means $s_{\pi_{*}}=o(\sqrt{n/log^{3}p})$, where
$s_{\pi_{*}}=\max\limits_{j \in G}{\|\pi_{*,j}\|_{0}}$; \\
(iv) $s_{\pi_{*}}\|\gamma_{*}\|_{0}+\|\theta_{*}\|_{0}=o(\sqrt{n}/ logp)$; \\
(v) there exist constants $\delta$ and $\kappa_{2}$ such that $E\vert u_{*}\vert^{2+\delta} \leq \kappa_{2}$; \\
(vi) there exist constants $K_{3} \ and \ K_{4} >0$, which only depend on the constants $\delta \  and\  \kappa_{2} $, satisfying that
 $
 \|\Sigma_{u_{*}}\gamma_{*}\|_{\infty}\geq \sqrt{n^{-1}logp}(K_{3}\|\gamma_{*}\|_{2}+K_{4})$.
\end{assumption}

Assumption \ref{assumption 3.1}(i) also considers  Gaussian design to simplify proofs.
Assumption \ref{assumption 3.1}(iii) shows that the  sparsity of the column of the matrix $\pi_{*}$.
Assumption \ref{assumption 3.1}(iv) imposes a sparse structure while the
non-sparseness  of model parameters $\beta_{A}$ and $\beta_{B}$ are allowed.
 For example,  $\beta_{A,G}=-\beta_{B,G}=\beta_{*}$ for a dense $\beta_{*}$ with $\|\beta_{*}\|_{0}=o(\sqrt{logp})$, and $\|\beta_{A,-G}\|_{0}=o(\sqrt{n}/logp)$, $\|\beta_{B,-G}\|_{0}=o(\sqrt{n}/logp)$.
Then $\beta_{A,G}+\beta_{B,G}=0$ and $\|\theta_{*}\|_{0}=0+\|\beta_{A,-G}\|_{0}+\|\beta_{B,-G}\|_{0}=o(\sqrt{n}/log p)$. Therefore, $s_{\pi_{*}}\|\gamma_{*}\|_{0}+\|\theta_{*}\|_{0}=o(\sqrt{n}/logp)$.
Assumption \ref{assumption 3.1}(vi) imposes a $l_{\infty}$ norm lower bound on the product of $\Sigma_{u_{*}}$ and $\gamma_{*}$.
As in the one-sample case, we give general Assumption \ref{assumption 3.1}(ii) and (v).
The following theorem shows that the test $S_{n}$ has a Type I error probability asymptotically equal to the nominal level $\alpha$.
\begin{theorem}\label{theorem 3.1}
 Let Assumption \ref{assumption 3.1}(i)-(iii) hold, when $n, p\rightarrow \infty$ with $logp=o(\sqrt{n})$, then under null hypothesis in \eqref{result hypothesis},
 \begin{equation}
  P(S_{n}>\Gamma^{-1}(1-\alpha;\hat{D}))\rightarrow\alpha, \forall\alpha\in(0,1)
 \tag{3.9},
 \end{equation}
  where $\Gamma^{-1}(1-\alpha;\hat{D})$  is the $1-\alpha$ quantile of the $\Gamma(x; \hat{D})$.
\end{theorem}
Similar to Theorem \ref{theorem 2.1}, Theorem \ref{theorem 3.1} does not make any model sparsity assumptions or require the error distribution.
Next, we consider the power property of  $S_{n}$.

 \begin{theorem} \label{theorem3.2}
Let Assumption \ref{assumption 3.1} holds, under the alternative hypothesis in \eqref{result hypothesis}, when $n, p\rightarrow \infty$ with $logp=o(\sqrt{n})$, then the test in Algorithm $\ref{2}$ is asymptotically  powerful, that is
 \begin{equation}
 P(S_{n}>\Gamma^{-1}(1-\alpha;\hat{D}))\rightarrow 1 , \notag \   {\forall} \alpha \in (0,1).
 \end{equation}
 \end{theorem}
 Theorem \ref{theorem3.2} shows that the test of group parameters  with  moment method in two samples is asymptotically powerful, that is, the probability of Type $\mathrm{II}$ error  tends to 0.
\section{Numerical results} \label{section 4}
In this section, we present finite-sample performance of the accuracy of the proposed methods.
In all simulations of one- and two-sample cases, we set $n=200$, $p=500$ and the nominal size of all the tests is 5\%. The rejection probabilities are based on 100 repetitions. For application purposes, we recommend choosing the tuning parameters as $\eta=0.5\sqrt{\frac{log p}{n}}$ and $\rho_{0}=0.01$. This is a general choice and we show in the simulations below that they provide good results.
\subsection{One-Sample Experiments }
Consider the model \eqref{1.1}, for simplicity of presentation,  we only test the following null hypothesis
\begin{equation}
 H_{0,G}: \beta_{G}=\beta_{G}^{0} \notag,
\end{equation}
 with $G=\{1,2,3\}$, $i.e.$ the first three components of $\beta$.
  We show the results for three different Gaussian designs as follows. \\
 \textit{Example 1.} Here we consider the standard Toeplitz  design where the rows of $X$ are drawn as an independent and identically distributed (i.i.d) random draws from a multivariate Gaussian distribution $N(0,\Sigma_{X})$, with covariance matrix $(\Sigma_{X})_{i,j}=0.4^{\vert i-j \vert}$. \\
\textit{Example 2.} In this example, we consider uncorrelated design where the rows of $X$ are i.i.d draws from  $N(0,\Sigma_{X})$, where $(\Sigma_{X})_{i,j}$ is 1 for $i=j$ and is 0 for $i\neq j$. \\
\textit{Example 3.} In this case, we consider a non-sparse design matrix with equal correlations among the features. Namely, the rows of $X$ are i.i.d draws from  $N(0,\Sigma_{X})$, where $(\Sigma_{X})_{i,j}$ is 1 for $i=j$ and is 0.4 for $i\neq j$.

 We also consider two specifications of the error distribution. \\
(a) In the light-tailed case, $\varepsilon$ is taken from a standard normal distribution.\\
(b) In the heavy-tailed case, the $\varepsilon$ is taken from Student's t-distribution with 3 degrees of freedom.

Let $s=\|\beta\|_{0}$ denotes the size of the model sparsity and we vary simulations setting from extremely sparse $s=3$ to extremely large $s=p$. For sparsity $s$, we set the model parameters as $\beta_{j}=\frac{3}{\sqrt{s}}$, $1\leq j \leq s$ and $\beta_{j}=0$, $j>s$.
To illustrate the non-sparse adaptability of MDS on group tests and the robustness of the moment method, we compare our method (i.e. the moment method combined with  MDS, MMDS) with the moment method combined with  debiased Lasso   (MDL) and the Wald test combined with  debiaed Lasso (WDL)  proposed by \cite{2014On}.
The implementation of MDL is as follows. We  compute the debiased estimators $\hat{\theta}_{debias}$ and $\hat{\pi}_{debias}$, so that $\hat{\sigma}_{e}=\|V-W\hat{\theta}_{debias}\|_{2}/\sqrt{n}$ and
$\hat{Q}_{debias}=\frac{1}{n}\sum\limits_{i=1}^{n}\hat{u}_{i}\hat{u}_{i}^{T}$, where $\hat{u}_{i}$ is the $i$-th row of  $\hat{U}=Z-W\hat{\pi}_{debias}$. Then the test is to reject $H_{0}$ if and only if
\begin{equation}
n^{-\frac{1}{2}}\hat{\sigma}_{e}^{-1}\|(Z-W\hat{\pi}_{debias})^{T}(V-W\hat{\theta}_{debias})\|_{\infty}> \Gamma^{-1}(1-\alpha; \hat{Q}_{debias}) \notag.
\end{equation}
The  WDL method is implemented as follows. We first compute the debiased estimator $\hat{\beta}_{debias}$ and  the nodewise Lasso estimator $\hat{\Theta}_{Lasso}$ for the precision matrix as in \citep{2014On}. Then the test is to reject $H_{0}$ if and only if
\begin{equation}
\max \limits_{j\in G}\frac{\sqrt{n}\vert \hat{\beta}_{debias;j}-\beta_{j}^{0}\vert}{\sigma_{\varepsilon}\sqrt{\hat{\Omega}_{j,j}}}
> \Psi^{-1}(1-\alpha) \notag
\end{equation}
where $\Psi=\max \limits_{j\in G}\frac{\vert W_{j} \vert}{\sigma_{\varepsilon}\sqrt{\hat{\Omega}_{j,j}}}$, with $W=\hat{\Theta}_{Lasso}X^{T}\varepsilon/\sqrt{n}\sim N_{n}(0, \sigma_{\varepsilon}^{2}\hat{\Omega})$, $\hat{\Omega}=\hat{\Theta}\hat{\Sigma}\hat{\Theta}^{T}$ and $\hat{\Sigma}=X^{T}X/n$.

We collect the results in Table \ref{Table 1}, where "light+Toeplitz" indicates light-tailed error distribution and Toeplitz design, others are similar. We can clearly see the instability of the methods of MDL and WDL, which means that the probabilities of the type I error are much higher than the nominal level $\alpha$  for both sparse and non-sparse models.
Conversely, when the sparsity of the model is equal to $s=p$, the type I error probability of MMDS remains stable. That is true  even if we change the correlation among the features and   error distribution.
However, the probabilities of type I error are very large for the MDL in the case of equal correlation of  feature  and for the WDL in the case of heavy-tailed error distribution, which indicates that these two methods are not robust to the covariance setting of the design matrix or error distribution.
\begin{table}[h]
\begin{center}
\begin{minipage}{\textwidth}
\caption{Type I error rate of MMDS , MDL and WDL in one sample}\label{Table 1}
\begin{tabular*}{\textwidth}{@{\extracolsep{\fill}}lccc|ccc|ccc@{\extracolsep{\fill}}}
\toprule%
& \multicolumn{3}{@{}c@{}}{Light+Toeplitz} & \multicolumn{3}{@{}c@{}}{Light+Noncorrelation} & \multicolumn{3}{@{}c@{}}{Light+Equal correlation} \\ \cmidrule{2-4}\cmidrule{5-7}\cmidrule{8-10}%
Method & MMDS       & MDL       & WDL       & MMDS           & MDL         & WDL        & MMDS            & MDL         & WDL \\
\midrule
s=3  & 0.04          & 0.31      & 0.28      & 0.05           & 0.30        & 0.08       & 0.04            & 0.99        &0.29 \\
s=5   & 0.03         & 0.38      & 0.35     & 0.03           & 0.28        & 0.16       & 0.04            & 1.00        & 0.28         \\
s=10  & 0.03         & 0.46      & 0.43     & 0.03           & 0.31        & 0.19       & 0.04            & 1.00         & 0.39        \\
s=20  & 0.04         & 0.45      & 0.51     & 0.04           & 0.26        & 0.50       & 0.03            & 1.00         & 0.48        \\
s=50  & 0.05         & 0.35      & 0.72     & 0.04           & 0.24        & 0.75       & 0.08            & 1.00         &0.73        \\
s=100 & 0.04         & 0.28      & 0.78     & 0.03           & 0.27        & 0.76       & 0.05            & 1.00        & 0.87      \\
s=n   & 0.03         & 0.28      & 0.69     & 0.07           & 0.24        & 0.71       & 0.05            & 1.00        &0.59        \\
s=p   & 0.03         & 0.35      & 0.67     & 0.03           & 0.23        & 0.57       & 0.05            & 1.00         & 0.67         \\
\hline
& \multicolumn{3}{@{}c@{}}{Heavy+Toeplitz} & \multicolumn{3}{@{}c@{}}{Heavy+Noncorrelation} & \multicolumn{3}{@{}c@{}}{Heavy+Equal correlation} \\ \cmidrule{2-4}\cmidrule{5-7}\cmidrule{8-10}%
Method & MMDS       & MDL       & WDL       & MMDS           & MDL         & WDL        & MMDS            & MDL         & WDL \\
\midrule
s=3   & 0.04         & 0.33      & 0.98     & 0.03           & 0.26        & 0.98       & 0.04            & 0.81         & 0.98         \\
s=5   & 0.05         & 0.31      & 0.99     & 0.05           & 0.16        & 0.98       & 0.03            & 0.81         & 0.98         \\
s=10  & 0.04         & 0.31      & 0.97     & 0.06           & 0.16        & 0.99       & 0.03            & 0.76         & 0.99        \\
s=20  & 0.04         & 0.32      & 1.00     & 0.04           & 0.26        & 0.99       & 0.06            & 0.71         & 1.00         \\
s=50  & 0.04         & 0.36      & 0.97     & 0.04           & 0.13        & 0.98       & 0.06            & 0.86         & 0.98        \\
s=100 & 0.04         & 0.30      & 0.94     & 0.06           & 0.19        & 0.93       & 0.06            & 0.92         & 0.96        \\
s=n   & 0.04         & 0.31      & 0.88     & 0.04           & 0.23        & 0.89       & 0.07            & 0.93         & 0.90        \\
s=p   & 0.04         & 0.38      & 0.73     & 0.03           & 0.30        & 0.76       & 0.08            & 0.92         & 0.84 \\
\botrule
\end{tabular*}
\end{minipage}
\end{center}
\end{table}

\begin{figure}[h]%
\centering
\includegraphics[width=0.9\textwidth]{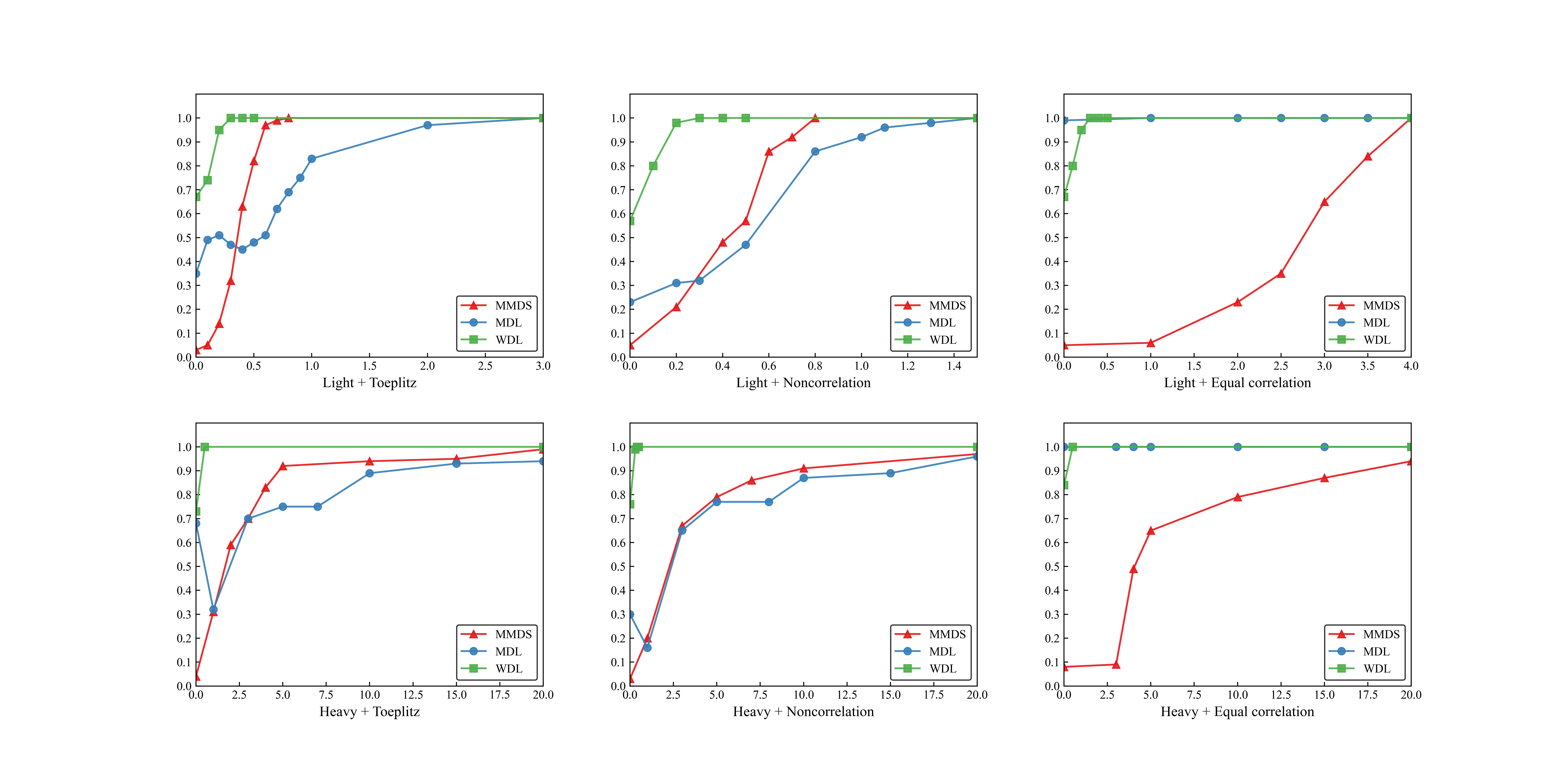}
\caption{Power curves of competing methods under different settings of design matrix and error distribution in one sample}\label{111}
\end{figure}

We also compare the power properties of MMDS with the MDL and WDL  in one sample.
The results are collected in Figure \ref{111}, which presents full power curves with various values of $h$, that measures the magnitude of deviations from the null hypothesis. Therefore, $h=0$ corresponds to Type I error  whereas other points on the curves correspond to Type II error ($h\neq 0$).
In all the power curves, we use a dense model with sparsity $p$, $i.e.$   $\|\beta\|_{0}=p$.
In figure \ref{111}, the first row corresponds to the light-tailed error distribution, and the second row corresponds to the heavy-tailed error distribution.
The two plots in the first column correspond to the Toeplitz designs, where we clearly observe that the MMDS outperforms both MDL and WDL by providing firm Type I error probability and also reaching full power quickly.
The second column corresponds to the uncorrelated design. The MDL and WDL  have  type I error probabilities over 0.05 whereas the MMDS still provides valid inference.
The last column corresponds to the design with equal correlation. The MDL and WDL  completely break down with Type I error probability being close to 1.
In conclusion, our method are stable in different settings of design matrices and model errors, while existing methods cannot control the probabilities of Type I error or Type II error.
\subsection{Two-Sample Experiments }
Consider the models \eqref{1.3} and \eqref{1.4} and errors $\varepsilon_{A}$ and $\varepsilon_{B}$.
The null hypothesis  is
\begin{equation}
H_{0}: \beta_{A,G}=\beta_{B,G}    \notag,
\end{equation}
with $G=\{1,2,3 \}$. Our goal is to test whether the first three components of the two model parameters are equal.
Let the design matrix $X_{A}$ be taken from Example 1-Example 3 introduced in one sample, and the covariance matrix of $X_{B}$ satisfies that $\Sigma_{B}=c\Sigma_{A}$. For simplicity of representation, we only consider the case of $c=2$. In addition, we also consider two specifications for the distributions of $\varepsilon_{A}$ and $\varepsilon_{B}$. \\
(a) In the light-tail case, $\varepsilon_{A}$ and $\varepsilon_{B}$ are drawn from the standard normal distribution. \\
(b) In the heavy-tail case, $\varepsilon_{A}$ and $\varepsilon_{B}$ are drawn from the Student's t-distribution  as described in one Sample.

Let $s_{A}=\|\beta_{A}\|_{0}$ denotes the size of the model sparsity and we vary simulations setting from extremely sparse $s_{A}=3$ to extremely large $s_{A}=p$. For sparsity $s_{A}$, we set $\beta_{A}$ as $\beta_{A,j}=\frac{3}{\sqrt{s_{A}}}$, $1\leq j \leq s_{A}$ and $\beta_{A,j}=0$, $j>s_{A}$. And let $\beta_{B}=\beta_{A}+(h,0,...,0)^{\top}$. The null hypothesis $H_{0}$ in \eqref{1.5} corresponds to $h=0$ and alternative hypothesis corresponds to $h\neq0$.

We compare our method (MMDS) with MDL and  WDL.  The implementation of MDL is as follows in two-sample case.  We first compute the debiased  estimators $\hat{\theta}_{*}$ and $\hat{\pi}_{*}$, so that $\hat{\sigma}_{\varepsilon_{*}}=\|Y-W\hat{\theta}_{*}\|_{2}/\sqrt{n}$ and
$\hat{D}_{debias}=\frac{1}{n}\sum\limits_{i=1}^{n}\hat{u}_{*,i}\hat{u}_{*,i}^{T}$, where $\hat{u}_{*,i}$ is the $i$-th row of  $\hat{U}_{*}=Z-W\hat{\pi}_{*}$. Then the test is to reject $H_{0}$ if and only if
\begin{align}
n^{-\frac{1}{2}}\hat{\sigma}_{\varepsilon_{*}}^{-1}\|(Z-W\hat{\pi}_{*})^{T}(Y-W\hat{\theta}_{*})\|_{\infty}> \Gamma^{-1}(1-\alpha; \hat{D}_{debias}) \notag.
\end{align}
The  WDL method in two-sample case is implemented as follows. We first compute the debiased estimators $\hat{\beta}_{debias;A}$ , $\hat{\beta}_{debias;B}$ and  the nodewise Lasso estimators $\hat{\Theta}_{Lasso;A}$ , $\hat{\Theta}_{Lasso;B}$ for the precision matrices as in \cite{2014On}. Then the test is to reject $H_{0}$ if and only if
\begin{align}
\max \limits_{j\in G}\frac{\sqrt{n} \vert \hat{\beta}_{debias;A,j}-\hat{\beta}_{debias;B,j} \vert}{\sigma_{\varepsilon,A}\sqrt{\hat{\Omega}_{A,j,j}}+\sigma_{\varepsilon,B}\sqrt{\hat{\Omega}_{B,j,j}}}
> \Phi^{-1}(1-\alpha) \notag
\end{align}
where $\Phi=\max \limits_{j\in G}\frac{\vert W_{A,j} - W_{B,j} \vert}{\sigma_{\varepsilon_{A}}\sqrt{\hat{\Omega}_{A,j,j}}+\sigma_{\varepsilon_{B}}\sqrt{\hat{\Omega}_{B,j,j}}}$, with $W_{A}=\hat{\Theta}_{Lasso;A}X_{A}^{T}\varepsilon_{A}/\sqrt{n}\sim N_{n}(0, \sigma_{\varepsilon_{A}}^{2}\hat{\Omega}_{A})$, $W_{B}=\hat{\Theta}_{Lasso;B}X_{B}^{T}\varepsilon_{B}/\sqrt{n}\sim N_{n}(0, \sigma_{\varepsilon_{B}}^{2}\hat{\Omega}_{B})$, $\hat{\Omega}_{A}=\hat{\Theta}_{Lasso;A}\hat{\Sigma}_{A}\hat{\Theta}_{Lasso;A}^{T}$, $\hat{\Omega}_{B}=\hat{\Theta}_{Lasso;B}\hat{\Sigma}_{B}\hat{\Theta}_{Lasso;B}^{T}$, $\hat{\Sigma}_{A}=X_{A}^{T}X_{A}/n$ and $\hat{\Sigma}_{B}=X_{B}^{T}X_{B}/n$.

The summary of the size results of two samples is presented in Table \ref{Table 2}, where we can see the non-robustness of MDL and WDL for the design matrices and error distributions, with the type I error probabilities being higher than the nominal level $\alpha$.
Conversely, even though the sparsity of the model is $p$ and the error has a heavy-tailed  distribution, the type I error probability of MMDS remains stable.
\begin{table}[h]
\begin{center}
\begin{minipage}{\textwidth}
\caption{Type I error rate of MMDS, MDL and WDL in two sample}\label{Table 2}
\begin{tabular*}{\textwidth}{@{\extracolsep{\fill}}lccc|ccc|ccc@{\extracolsep{\fill}}}
\toprule%
& \multicolumn{3}{@{}c@{}}{Light+Toeplitz} & \multicolumn{3}{@{}c@{}}{Light+Noncorrelation} & \multicolumn{3}{@{}c@{}}{Light+Equal correlation} \\ \cmidrule{2-4}\cmidrule{5-7}\cmidrule{8-10}%
Method & MMDS       & MDL       & WDL       & MMDS           & MDL         & WDL        & MMDS            & MDL         & WDL \\
\midrule
      & MMDS      & MDL      & WDL     & MMDS        & MDL        & WDL        & MMDS         & MDL         & WDL        \\
s=3   & 0.04         & 0.15      & 0.03     & 0.04           & 0.14        & 0.01       & 0.07            & 0.54         & 0.07         \\
s=5   & 0.02         & 0.24      & 0.04     & 0.04           & 0.12        & 0.04       & 0.02            & 0.71         & 0.09        \\
s=10  & 0.06         & 0.12      & 0.08     & 0.03           & 0.18        & 0.11       & 0.03            & 0.84         & 0.16        \\
s=20  & 0.03         & 0.11      & 0.10     & 0.08           & 0.10        & 0.23       & 0.04            & 0.87         & 0.23        \\
s=50  & 0.04         & 0.13      & 0.51     & 0.04           & 0.10        & 0.50       & 0.02            & 0.89         & 0.41        \\
s=100 & 0.03         & 0.12      & 0.60     & 0.03           & 0.17        & 0.52       & 0.03            & 0.91         & 0.55        \\
s=n   & 0.03         & 0.08      & 0.61     & 0.08           & 0.15        & 0.56       & 0.02            & 0.90         & 0.63        \\
s=p   & 0.04         & 0.16      & 0.50     & 0.06           & 0.07        & 0.41       & 0.03            & 0.95         & 0.57       \\
\hline
& \multicolumn{3}{@{}c@{}}{Heavy+Toeplitz} & \multicolumn{3}{@{}c@{}}{Heavy+Noncorrelation} & \multicolumn{3}{@{}c@{}}{Heavy+Equal correlation} \\ \cmidrule{2-4}\cmidrule{5-7}\cmidrule{8-10}%
Method & MMDS       & MDL       & WDL       & MMDS           & MDL         & WDL        & MMDS            & MDL         & WDL \\
\midrule
s=3   & 0.06         & 0.07      & 0.09     & 0.03           &   0.08      & 0.02       & 0.03            & 0.99         & 0.03        \\
s=5   & 0.06         & 0.07      & 0.10     & 0.03           &   0.08      & 0.03       & 0.07            & 0.40         & 0.08        \\
s=10  & 0.06         & 0.05      & 0.10     & 0.06           &   0.08      & 0.03       & 0.08            & 0.47         & 0.17       \\
s=20  & 0.05         & 0.07      & 0.13     & 0.04           &   0.06      & 0.11       & 0.02            & 0.47         & 0.10      \\
s=50  & 0.07         & 0.06      & 0.26     & 0.04           &   0.08      & 0.10       & 0.03            & 0.50         & 0.10        \\
s=100 & 0.06         & 0.10      & 0.22     & 0.04           &    0.07     & 0.09       & 0.07            & 0.50         & 0.17       \\
s=n   & 0.04         & 0.15      & 0.27     & 0.04           &    0.08     & 0.15       & 0.03            &  0.54        & 0.16        \\
s=p   & 0.05         & 0.10      & 0.16     & 0.06           &    0.07     & 0.10       & 0.05            & 0.70         & 0.26    \\
\botrule
\end{tabular*}
\end{minipage}
\end{center}
\end{table}

\begin{figure}[h] \footnotesize
\centering
\includegraphics[width=0.9\textwidth]{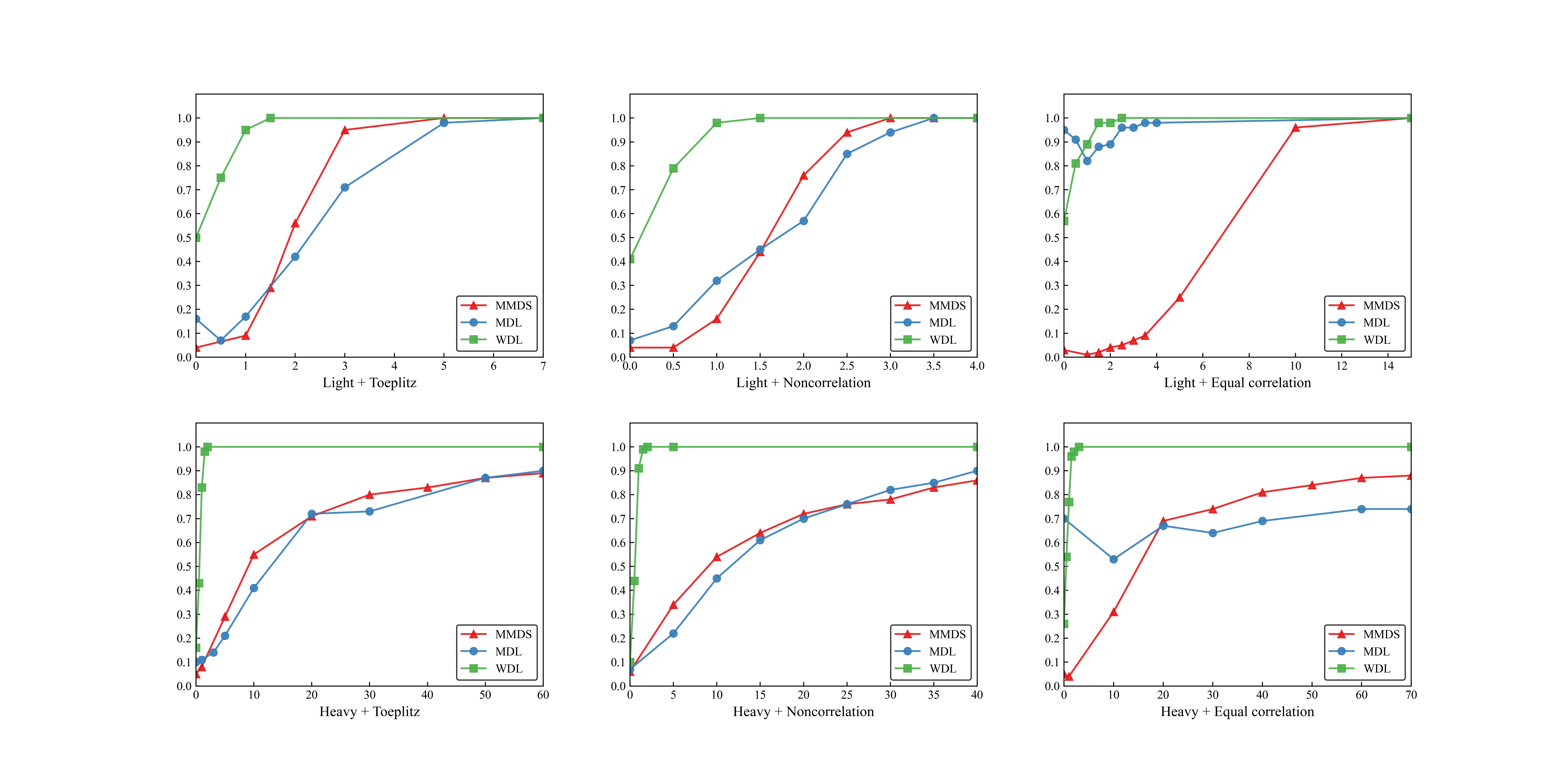}
\caption{Power curves of competing methods under different settings of design matrix and error distribution in two samples}\label{222}
\end{figure}

We also compare the power properties of MMDS with the MDL and WDL  in two samples.
In all the power curves of two samples, we use  dense models with sparsity $p$, $i.e.$   $\|\beta_{A}\|_{0}=\|\beta_{B}\|_{0}=p$.
The results are collected in Figure \ref{222}.  It can be seen that the speed of reaching full power of MMDS is comparable to that of MDL in most settings. However, MDL fails to control the probability of  Type I error.
In addition, WDL reaches full power the fastest, but its Type I error is well above the nominal level $\alpha$. Especially in the case of a light-tailed error distribution, the probability of type I error probability is around 0.5.

\section{Conclusion}
In this paper, we consider simultaneous test of  group  parameters under non-sparse high-dimensional linear model in one- and two-sample cases. We employ reconstruction and convolutional regression methods to transform the null hypothesis into a testable moment condition, and estimate the unknown parameters in it with   Modified Dantzig Selector.  By plugging-in the estimators, we construct the test statistic.
In addition, we demonstrate the excellent properties of the test in controlling Type I and Type II errors through theoretical and simulation experiments.
The above procedure does not make model sparsity assumptions, which is difficult to verify and satisfy, so it can solve more practical problems.
\backmatter

\bmhead{Supplementary information}
Supplement to "Simultaneous Inference  in Non-Sparse High-Dimensional Linear Models".
The detailed proofs about the asymptotic distributions of test statistics constructed in one- and two-sample models are given. In addition,  we also give detailed proofs of the power properties of the tests. Technical lemmas are also proved in the supplement.
\bmhead{Acknowledgments}
This paper is supported by  National Social Science Fund project of China (21BTJ045).





\bibliography{sn-bibliography}

\end{document}